\documentclass[12pt]{amsart}
\usepackage{amscd,amsmath,amssymb,amsthm,enumerate,times,ifthen}
\usepackage{epsfig,epic,eepic}
\usepackage{graphicx}
\usepackage[utf8]{inputenc}
\usepackage[T1]{fontenc}
\usepackage{tikz}
\usetikzlibrary{matrix,arrows}
\usepackage[mathscr]{eucal}
\newtheorem{theorem}{Theorem}
\newtheorem{Theorem}{Theorem}

\def\Thm#1#2{\begin{theorem}\label{T#1}#2\end{theorem}}

\def\thm#1{Theorem~\ref{T#1}}

\newtheorem{Lemma}{Lemma}
\def\Lem#1#2{\begin{Lemma}\label{L#1}#2\end{Lemma}}
\def\lem#1{Lemma~\ref{L#1}}

\newtheorem{Prposition}{Proposition}

\newtheorem{Corollary}{Corollary}
\def\Cor#1#2{\begin{Corollary}\label{C#1}#2\end{Corollary}}
\def\cor#1{Corollary~\ref{C#1}}

\newtheorem*{theoremn}{Theorem}
\def\Thmn#1{\begin{theoremn}#1\end{theoremn}}

\newtheorem*{lemman}{Lemma}

\newtheorem*{corn}{Corollary}

\theoremstyle{definition}
\newtheorem*{Def}{Definition}

\newcommand{\N}{\mathbb N}

\newcommand{\R}{\mathbb R}

\newcommand{\al}{\alpha}
\newcommand{\be}{\beta}

\newcommand{\de}{\delta}

\newcommand{\eps}{{\varepsilon}}
\newcommand{\diam}{\mathop{\hbox{\rm diam}}}
\renewcommand{\phi}{\varphi}
\newcommand{\eq}[1]{\eq{#1}}
\newcommand{\Eq}[2]{\ifthenelse{\equal{#1}{*}}
{\begin{equation*}\begin{aligned}#2\end{aligned}\end{equation*}}
{\begin{equation}\label{#1}\begin{aligned}#2\end{aligned}\end{equation}}}

\begin{document}

\date{\today}

%\begin{flushright}
%\emph{Journal of Fixed Point Theory and Applications}
%\end{flushright}

\title{A contraction principle in semimetric spaces}

\author[M. Bessenyei]{Mih\'aly Bessenyei}
\author[Zs. P\'ales]{Zsolt P\'ales}

\address{Institute of Mathematics,
University of Debrecen, H-4010 Debrecen, Pf.\ 12, Hungary}

\email{besse@science.unideb.hu}
\email{pales@science.unideb.hu}



\thanks{This research was realized in the frames of T\'AMOP 4.2.4. A/2-11-1-2012-0001 
``National Excellence Program Elaborating and operating an inland student and researcher personal
support system''. The project was subsidized by the European Union and co-financed by the European 
Social Fund. This research was also supported by the Hungarian Scientific 
Research Fund (OTKA) Grant NK 81402.}

\begin{abstract}
A branch of generalizations of the Banach Fixed Point Theorem replaces contractivity by a weaker but still
effective property. The aim of the present note is to extend the contraction principle in this spirit
for such complete semimetric spaces that fulfill an extra regularity property. The stability of fixed 
points is also investigated in this setting. As applications, fixed point results are presented for several 
important generalizations of metric spaces.
\end{abstract}

\maketitle

\section{Introduction}

Although the contraction principle appears partly in the method of successive approximation in the works of
Cauchy \cite{Cau1835}, Liouville \cite{Lio1837}, and Picard \cite{Pic1890}, its abstract and powerful version
is due to Banach \cite{Ban22}. This form of the contraction principle, commonly quoted as the Banach Fixed Point
Theorem, states that any contraction of a complete metric space has exactly one fixed point. Until now, this
seminal result has been generalized in several ways and some of these generalizations initiated new branches
in the field of Iterative Fixed Point Theory. The books by Granas and Dugundji \cite{GraDug03}, and by Zeidler
\cite{Zei86} give an excellent and detailed overview of the topic.

Some generalizations of Banach's fundamental result replace contractivity by a weaker but still effective
property; for example, the self-mapping $T$ of a metric space $X$ is supposed to satisfy
\Eq{main}{
 d(Tx,Ty)\le\phi\bigl(d(x,y)\bigr)\qquad(x,y\in X).}

To the best of our knowledge, the assumption above appears first in the paper by Browder \cite{Brow68} and
by Boyd and Wong \cite{BoyWon69}. One of the most important result in this setting was obtained by Matkowski
\cite{Mat75} who established the following statement.

\Thmn{Assume that $(X,d)$ is a complete metric space and $\phi\colon\R_+\to\R_+$ is a monotone increasing 
function
such that the sequence of iterates $(\phi^n)$ tends to zero pointwise on the positive half-line. If $T\colon 
X\to X$
is a mapping satisfying \eqref{main}, then it has a unique fixed point in $X$.}

Another branch of generalizations of Banach's principle is based on relaxing the axioms of the metric space. 
For an account of such developments, see the monographs by Rus \cite{Rus01}, by Rus, Petru\c{s}el and Petru\c{s}el 
\cite{RusPetPet08}, and by Berinde \cite{Ber07}. The aim of this note is to extend the contraction principle
combining these two directions: The main result is formulated in the spirit of Matkowski, while the underlying
space is a complete semimetric space that fulfills an extra regularity condition. These kind of spaces involve
standard metric, ultrametric and inframetric spaces. The stability of fixed points is also investigated in this
general setting.

\section{Conventions and Basic Notions}

Throughout this note, $\R_+$ and $\overline{\R}_+$ stand for the set of all nonnegative and extended nonnegative
reals, respectively. The \emph{iterates} of a mapping $T\colon X\to X$ are defined inductively by the 
recursion
$T^1=T$ and $T^{n+1}=T\circ T^n$.

Dropping the third axiom of Fr\'echet \cite{Fre1906}, we arrive at the notion of semimetric spaces: Under a
\emph{semimetric space} we mean a pair $(X,d)$, where $X$ is a nonempty set, and $d\colon X\times X\to\R_+$
is a nonnegative and symmetric function which vanishes exactly on the diagonal of the Cartesian product
$X\times X$. In semimetric spaces, the notions of \emph{convergent} and \emph{Cauchy sequences}, as like as
(open) \emph{balls} with given center and radius, can be introduced in the usual way. For an open ball with
center $p$ and radius $r$, we use the notation $B(p,r)$. Under the \emph{diameter} of $B(p,r)$ we mean the
supremum of distances taken over the pairs of points of the ball. Under the \emph{topology} of a semimetric
space we mean the topology induced by open balls.

\begin{Def}
Consider a semimetric space $(X,d)$. We say that $\Phi\colon\overline{\R}_+^2\to\overline{\R}_+$ is a 
\emph{triangle function} for $d$, if $\Phi$ is symmetric and monotone increasing in both of its arguments, 
satisfies $\Phi(0,0)=0$ and, for all $x,y,z\in X$, the generalized triangle inequality
\Eq{*}{
 d(x,y)\le\Phi\bigl(d(x,z),d(y,z)\bigr)}
holds.
\end{Def}

The construction below plays a key role in the further investigations. For a semimetric space $(X,d)$, define
the function
\Eq{basic}{
 \Phi_d(u,v):=\sup\{d(x,y)\mid\exists p\in X: d(p,x)\le u,\,d(p,y)\le v\}\qquad(u,v\in\overline{\R}_+).}
Simple and direct calculations show, that $\Phi_d\colon\overline{\R}^2_+\to\overline{\R}_+$ is a triangle 
function for $d$. This function is called the \emph{basic triangle function}. Note also, that the basic 
triangle function is optimal in the following sense: If $\Phi$ is a triangle function for $d$, then 
$\Phi_d\le\Phi$ holds.

Obviously, metric spaces are semimetric spaces with triangle function $\Phi(u,v):=u+v$. Ultrametric spaces are
also semimetric spaces if we choose $\Phi(u,v):=\max\{u,v\}$. Not claiming completeness, we present here
some further examples that can be interpreted in this framework:
\begin{itemize}
 \item $\Phi(u,v)=c(u+v)$ ($c$-relaxed triangle inequality);
 \item $\Phi(u,v)=c\max\{u,v\}$ ($c$-inframetric inequality);
 \item $\Phi(u,v)=(u^p+v^p)^{1/p}$ ($p$th-order triangle inequality, where $p>0$).
\end{itemize}

Briefly, each semimetric space $(X,d)$ can be equipped with an optimal triangle function attached to $d$; this
triangle function provides an inequality that corresponds to and plays the role of the classical triangle
inequality. In this sense, semimetric spaces are closer relatives of metric spaces then the system of axioms
suggests.

Throughout the present note, we shall restrict our attention only to those semimetric spaces whose basic
triangle function is continuous at the origin. These spaces are termed \emph{regular}. Clearly,
the basic triangle function of a regular semimetric space is bounded in a neighborhood of the origin. The
importance of regular semimetric spaces is enlightened by the next important technical result.  

\Lem{regsem}{The topology of a regular semimetric space is Hausdorff. A convergent sequence in a regular
semimetric space has a unique limit and possesses the Cauchy property. Moreover, a semimetric space $(X,d)$
is regular if and only if
\Eq{diameter}{
 \lim_{r\to 0}\sup_{p\in X}\diam B(p,r)=0.}}

\begin{proof}
For the first statement, assume to the contrary that there exist distinct points $x,y\in X$ such that, for all
$r>0$, the balls $B(x,r)$ and $B(y,r)$ are not disjoint. The continuity and the separate monotonicity of the
basic triangle function guarantee the existence of $\de>0$ such that, for all $r<\de$, we have $\Phi_d(r,r)<d(x,y)$.
Therefore, if $p\in B(x,r)\cap B(y,r)$, we get the contradiction
\Eq{*}{
 d(x,y)\le\Phi_d\bigl(d(p,x),d(p,y)\bigr)\le\Phi_d(r,r)<d(x,y).}

The Hausdorff property immediately implies that the limit of a convergent sequence is unique. Assume now that
$(x_n)$ is convergent and tends to $x\in X$. Then the generalized triangle inequality implies the estimation
$d(x_n,y_m)\le\Phi_d\bigl(d(x,x_n),d(x,x_m)\bigr)$. The regularity of the underlying space provides that the
right-hand side tends to zero if we take the limit $n\to\infty$, yielding the Cauchy property.

If $(X,d)$ is a regular semimetric space, then the basic triangle function $\Phi_d$ is continuous at the origin.
Hence, for $\eps>0$, there exists $u_0,v_0>0$ such that $\Phi_d(u,v)<\eps$ whenever $0<u<u_0$ and $0<v<v_0$. Let
$r_0=\min\{u_0,v_0\}$. Fix $0<r<r_0$ and $p\in X$. If $x,y\in B(p,r)$, then, using the separate monotonicity of
triangle functions,
\Eq{*}{
 d(x,y)\le\Phi_d\bigl(d(p,x),d(p,y)\bigr)\le\Phi_d(r,r)<\eps}
follows. That is, $\diam B(p,r)\le\eps$ holds for all $p\in X$. Since $\eps$ is an arbitrary positive number,
we arrive at the desired limit property.

Assume conversely that the diameters of balls with small radius are uniformly small, and take sequences of
positive numbers $(u_n)$ and $(v_n)$ tending to zero. For fixed $n\in\N$, define $r_n=\max\{u_n,v_n\}$ and
take elements $p,x,y\in X$ satisfying $d(p,x)\le u_n$ and $d(p,y)\le v_n$. Then, $x,y\in B(p,r_n)$; therefore
\Eq{*}{
 d(x,y)\le\diam B(p,r_n)\le\sup_{p\in X}\diam B(p,r_n).}
Taking into account the definition of the basic triangle function $\Phi_d$ and the choice of $p,x,y$, we arrive
at the inequality
\Eq{*}{
 \Phi_d(u_n,v_n)\le\sup_{p\in X}\diam B(p,r_n).}
Here the right-hand side tends to zero as $n\to\infty$ by hypothesis, resulting the continuity of $\Phi_d$ at
the origin. 
\end{proof}

As usual, a semimetric space is termed to be \emph{complete}, if each Cauchy sequence of the space is convergent.
In view of the previous lemma, convergence and Cauchy property cannot be distinguished in complete and regular
semimetric spaces. 

In order to construct a large class of complete, regular semimetric spaces, we introduce the notion of 
equivalence of semimetrics. Given two semimetrics $d_1$ and $d_2$ on $X$, an increasing function 
$L\colon\overline{\R}_+\to\overline{\R}_+$ such that $L(0)=0$ and
\Eq{*}{
  d_1(x,y)\leq L(d_2(x,y)) \qquad(x,y\in X)
}
is called a \emph{Lipschitz modulus with respect to the pair $(d_1,d_2)$}. It is immediate to see that
the function $L_{d_1,d_2}\colon\overline{\R}_+\to\overline{\R}_+$ defined by
\Eq{*}{
  L_{d_1,d_2}(t):=\sup\{d_1(x,y)\mid x,y\in X,\,d_2(x,y)\leq t\} \qquad(t\in \overline{\R}_+)
}
is the smallest Lipschitz modulus with respect to $(d_1,d_2)$. The semimetrics $d_1$ and $d_2$ are said to be 
\emph{equivalent} if 
\Eq{*}{
  \lim_{t\to0+} L_{d_1,d_2}(t)=0\qquad\mbox{and}\qquad \lim_{t\to0+} L_{d_2,d_1}(t)=0
}
i.e., if $L_{d_1,d_2}$ and $L_{d_2,d_1}$ are continuous at zero.
It is easy to verify that, indeed, this notion of equivalence is an equivalence relation. For the proof of 
the transitivity one should use the inequality $ L_{d_1,d_3}\leq L_{d_1,d_2}\circ L_{d_2,d_3}$.

Our next result establishes that the convergence, completeness and the regularity of a semimetric space is 
invariant with respect to the equivalence of the semimetrics.

\Lem{equiv}{If $d_1$ and $d_2$ are semimetrics on $X$, then
\Eq{Lip}{
   \Phi_{d_1}\leq L_{d_1,d_2}\circ\Phi_{d_2}\circ(L_{d_2,d_1},L_{d_2,d_1}).
}
Provided that $d_1$ and $d_2$ are equivalent semimetrics, we have that 
\begin{enumerate}[(i)]
 \item a sequence converges to point in $(X,d_1)$ if and only if it converges to the 
same point in $(X,d_2)$;
 \item a sequence is Cauchy in $(X,d_1)$ if and only if it is Cauchy in $(X,d_2)$;
 \item $(X,d_1)$ is complete if and only if $(X,d_2)$ is complete;
 \item $(X,d_1)$ is regular if and only if $(X,d_2)$ is regular.
\end{enumerate}}

\begin{proof} Using the monotonicity properties, for $x,y,z\in X$, we have
\Eq{*}{
  d_1(x,y)\leq L_{d_1,d_2}(d_2(x,y)) 
          &\leq L_{d_1,d_2}\Big(\Phi_{d_2}(d_2(x,z),d_2(z,y))\big)\Big) \\
          &\leq L_{d_1,d_2}\Big(\Phi_{d_2}\big(L_{d_2,d_1}(d_1(x,z)),L_{d_2,d_1}(d_1(z,y))\big)\Big),
}
whence it follows that the map $L_{d_1,d_2}\circ\Phi_{d_2}\circ(L_{d_2,d_1},L_{d_2,d_1})$ is a triangle 
function for $d_1$, hence \eqref{Lip} must be valid.

For (i), let $(x_n)$ be a sequence converging to $x$ in $(X,d_1)$. Then $(d_1(x_n,x))$ is a null-sequence.
By the continuity of $L_{d_2,d_1}$ at zero, the right-hand side of the inequality 
\Eq{*}{
d_2(x_n,x)\leq L_{d_2,d_1}(d_1(x_n,x))} 
also tends to zero, hence $(d_2(x_n,x))$ is also a null-sequence. The reversed implication holds analogously.

The proof for the equivalence of the Cauchy property is completely similar.

To prove (iii), assume that $(X,d_1)$ is complete and let $(x_n)$ be a Cauchy sequence in $(X,d_2)$. Then,
by (ii), $(x_n)$ is Cauchy in $(X,d_1)$. Hence, there exists an $x\in X$ such that $(x_n)$ converges to $x$ 
in $(X,d_1)$. Therefore, by (i), $(x_n)$ converges to $x$ in $(X,d_2)$. This proves the completeness of 
$(X,d_2)$. The reversed implication can be verified analogously.

Finally, assume that $(X,d_2)$ is a regular semimetric space, which means that $\Phi_{d_2}$ is continuous at 
$(0,0)$. Using inequality \eqref{Lip}, it follows that $\Phi_{d_1}$ is also continuous at $(0,0)$ yielding 
that $(X,d_1)$ is regular, too.
\end{proof}

By the above result, if a semimetric is equivalent to a complete metric, then it is regular and complete.
With a similar argument that was followed in the above proofs, one can easily verify that equivalent 
semimetrics on $X$ generate the same topology.

\medskip

In the sequel, we shall need a concept that extends the notion of classical contractions to nonlinear ones.
This extension is formulated applying comparison functions fulfilling the assumptions of Matkowski.   

\begin{Def}
Under a \emph{comparison function} we mean a monotone increasing function $\phi\colon\R_+\to\R_+$ such
that the limit property $\lim_{n\to\infty}\phi^n(t)=0$ holds for all $t\ge 0$. Given a semimetric space
$(X,d)$ and a comparison function $\phi$, a mapping $T\colon X\to X$ is said to be \emph{$\phi$-contractive}
or a \emph{$\phi$-contraction} if it fulfills \eqref{main}.
\end{Def}

The statement of the next lemma is well-known, we provide its proof for the convenience of the reader. 

\Lem{mainprop}{If $\phi$ is a comparison function, then $\phi(t)<t$ for all positive $t$. If $(X,d)$ is a 
semimetric space and $T\colon X\to X$ is a $\phi$-contraction, then $T$ has at most one fixed point.}

\begin{proof}
For the first statement, assume at the contrary that $t\le\phi(t)$ for some $t>0$. Whence, by monotonicity and
using induction, we arrive at
\Eq{*}{
 t\le\phi(t)\le\phi^2(t)\le\cdots\le\phi^n(t).}
Upon taking the limit $n\to\infty$, the right hand-side tends to zero, contradicting to the positivity of $t$.
In view of this property, the second statement follows immediately. Indeed, if $x_0,y_0$ were distinct fixed
points of a $\phi$-contraction $T$, then we would arrive at 
\Eq{*}{
 t=d(x_0,y_0)=d(Tx_0,Ty_0)\le\phi(d(x_0,y_0))=\phi(t)<t.}
This contradiction implies $x_0=y_0$.
\end{proof}

\section{The Main Results}

The main results of this note is presented in two theorems. The first one is an extension of the Matkowski Fixed
Point Theorem \cite{Mat75} for complete and regular semimetric spaces. The most important ingredient of the proof
is that a domain invariance property remains true for sufficiently large iterates of $\phi$-contractions.  

\Thm{mainfix}{If $(X,d)$ is a complete regular semimetric space and $\phi$ is a comparison function, then every
$\phi$-contraction has a unique fixed point.}

\begin{proof}
Let $T\colon X\to X$ be a $\phi$-contraction and let $p\in X$ be fixed arbitrarily. Define the sequence 
$(x_n)$ by the
standard way $x_n:=T^np$. Observe that, for all fixed $k\in\N$, the sequence $\bigl(d(x_n,x_{n+k})\bigr)$ tends
to zero. Indeed, by the asymptotic property of comparison functions,
\Eq{*}{
 d(x_n,x_{n+k})=d(Tx_{n-1},Tx_{n+k-1})\le
  \phi\bigl(d(x_{n-1},x_{n+k-1})\bigr)\le\cdots\le
   \phi^n\bigl(d(p,T^kp)\bigr)\longrightarrow 0.}

We are going to prove that $(x_n)$ is a Cauchy sequence. Fix $\eps>0$. The continuity of the basic triangle 
function guarantees the existence of a neighborhood $U$ of the origin such that, for all $(u,v)\in U$, we 
have the inequality $\Phi_d(u,v)<\eps$. Or equivalently, applying the separate monotonicity, there exists 
some $\de(\eps)>0$ such that $\Phi_d(u,v)<\eps$ holds if $0\le u,v<\de(\eps)$. The asymptotic property of 
comparison functions allows us to fix an index $n(\eps)\in\N$ such that $\phi^{n(\eps)}(\eps)<\de(\eps)$ 
hold. Then, $\psi:=\phi^{n(\eps)}$ is a comparison function. Hence, if $0\le u,v<\min\{\eps,\de(\eps)\}$, 
\Eq{*}{
 \Phi_d\bigl(u,\psi(v)\bigr)\le\Phi_d(u,v)<\eps.}
Immediate calculations show, that the mapping $S:=T^{n(\eps)}$ is a $\psi$-contraction. 
Let the nonnegative integer $k$ and the points $x,y\in X$ be arbitrary. Then, applying the monotonicity 
properties of comparison functions and
their iterates,
\Eq{*}{
 d(T^kSx,T^kSy)\le\psi\bigl(d(T^kx,T^ky)\bigr)
  \le\psi\circ\phi^k\bigl(d(x,y)\bigr)
   \le\psi\bigl(d(x,y)\bigr)}
follows. This inequality immediately implies that $T^kS$ maps the ball $B(x,\eps)$ into itself if it makes small
perturbation on the center. Indeed, if $y\in B(x,\eps)$ and $d(x,T^kSx)<\de(\eps)$, the choices of $\eps$ and
$\de(\eps)$ moreover the separate monotonicity of the basic triangle function yield
\Eq{*}{
 d(x,T^kSy)
  \le\Phi_d\bigl(d(x,T^kSx),d(T^kSx,T^kSy)\bigr)
   \le\Phi_d\bigl(d(x,T^kSx),\psi(d(x,y)\bigr)<\eps.}

The properties of $(x_n)$ established at the beginning of the proof ensure that, for all nonnegative $k$, there
exists some $n_k\in\N$ such that the inequalities $d(x_n,T^kSx_n)<\de(\eps)$ hold whenever $n\ge n_k$. Choose
\Eq{*}{
 n_0=\max\{n_k\mid k=1,\ldots,n(\eps)\}.}
Then, taking into account the previous step, $T^kS\colon B(x_{n_0},\eps)\to B(x_{n_0},\eps)$ for 
$k=1,\ldots,n(\eps)$.
In particular, each iterates of $S$ is a self-mapping of the ball $B(x_{n_0},\eps)$. Let $n>n_0$ be an arbitrarily
given natural number. Then $n=mn(\eps)+k$, where $m\in\N$ and $k\in\{1,\ldots,n(\eps)\}$; hence, the definition of
$S$ leads to
\Eq{*}{
 T^nS=T^{mn(\eps)+k}S=T^kT^{mn(\eps)}S=T^kS^mS=T^kS^{m+1}.}
Therefore,
\Eq{*}{
 T^nS\bigl(B(x_{n_0},\eps)\bigr)
  =T^kS^{m+1}\bigl(B(x_{n_0},\eps)\bigr)
   \subset T^kS\bigl(B(x_{n_0},\eps)\bigr)\subset B(x_{n_0},\eps).} 
In other words, due to property \eqref{diameter} of \lem{regsem}, the sequence $(T^nSx_{n_0})$ is Cauchy and hence
so is $(x_n)$. The completeness implies, that it tends to some element $x_0$ of $X$. Our claim is that $x_0$ is a
fixed point of $T$. Applying the generalized triangle inequality,
\Eq{*}{
 d(x_0,Tx_0)&\le\Phi_d\bigl(d(x_0,x_{n+1}),d(Tx_0,x_{n+1})\bigr)\\
            &=\Phi_d\bigl(d(x_0,x_{n+1}),d(Tx_0,Tx_n)\bigr)\\
            &\le\Phi_d\bigl(d(x_0,x_{n+1}),\phi(d(x_0,x_n))\bigr)\\
            &\le\Phi_d\bigl(d(x_0,x_{n+1}),d(x_0,x_n)\bigr).}
Therefore,
\Eq{*}{
 d(x_0,Tx_0)\le\lim_{n\to\infty}\Phi_d\bigl(d(x_0,x_{n+1}),d(x_0,x_n)\bigr)=\Phi_d(0,0)=0}
follows. That is, $Tx_0=x_0$, as it was desired. To complete the proof, recall that a $\phi$-contraction 
may have at most one fixed point.
\end{proof}

Our second main result is based on \thm{mainfix} whose proof requires some additional ideas. It asserts
the stability of fixed points of iterates of $\phi$-contractions. 

\Thm{stabfix}{Let $(X,d)$ be a complete and regular semimetric space. If $(T_n)$ is a sequence of 
$\phi$-contractions converging pointwise to a $\phi$-contraction $T_0\colon X\to X$, then 
the sequence of the fixed points of $(T_n)$ converges to the unique fixed point of $T_0$.}

\begin{proof}
First we show that, for all $k\in\N$, $(T_n^k)$ converges to $T_0^k$ pointwise.
We prove by induction on $k$. By the assumption of the theorem, we have the statement for $k=1$.
Now assume that $T_n^k \to T_0^k$ pointwise. Then, for every $x\in X$, we have
\Eq{*}{
  d(T_n^{k+1}x,T_0^{k+1}x)
    &\leq \Phi_d\bigl(d(T_n^{k+1}x,T_n^{k}T_0x),d(T_n^{k}T_0x,T_0^{k+1}x)\bigr)\\
    &\leq \Phi_d\bigl(\phi^k(d(T_nx,T_0x)),d(T_n^{k}T_0x,T_0^{k}T_0x)\bigr)\\
    &\leq \Phi_d\bigl(d(T_nx,T_0x),d(T_n^{k}T_0x,T_0^{k}T_0x)\bigr)\to0
}
since $T_n\to T_0$, $T_n^k \to T_0^k$ pointwise and the semimetric $d$ is regular. This proves that
$(T_n^{k+1})$ converges to $T_0^{k+1}$ pointwise. 

The previous theorem ensures that, for all $n\in\N$, there exists
a unique fixed point $x_n$ of $T_n$ as well as a unique fixed point $x_0$ of $T_0$. 
To verify the statement of the theorem, assume indirectly that
\Eq{*}{
 \eps:=\limsup_{n\to\infty}d(x_0,x_n)>0.}
Choose $\delta>0$ such that $\Phi_d\bigl(\de,\de\bigr)<\eps$ hold, and then choose $m\in\N$ such 
that $\phi^m(2\eps)<\de$, finally denote $\psi:=\phi^m$. Using induction, one can easily prove that 
$S_n=T^m_n$ is a $\psi$-contraction and $x_n$ is a fixed point of $S_n$ for all $n\geq 0$. Furthermore,
in view of the previous step, the sequence $(S_nx_0)$ converges to $S_0x_0=x_0$. Hence, for large $n$,
$d(x_0,S_nx_0)<\delta$ and $d(x_n,x_0)<2\eps$. Therefore, for large $n$, we have
\Eq{*}{
 d(x_0,x_n)&\le\Phi_d\bigl(d(x_0,S_nx_0),d(S_nx_0,x_n)\bigr)\\
           &=\Phi_d\bigl(d(x_0,S_nx_0),d(S_nx_0,S_nx_n)\bigr)\\
           &\le\Phi_d\bigl(d(x_0,S_nx_0),\psi(d(x_n,x_0))\bigr)
            \le\Phi_d\bigl(\de,\psi(2\eps)\bigr).
}
Taking the limes superior, 
$\eps\leq\Phi_d\bigl(\de,\psi(2\eps)\bigr) \leq\Phi_d\bigl(\de,\de\bigr)<\eps$ follows, which is a 
contradiction. That is, the sequence of fixed points of $T_n$ tends to the fixed point of $T_0$, as it was 
stated.  
\end{proof}

If the semimetric $d$ is \emph{self-continuous}, that is, $d(x_n,y_n)\to d(x,y)$ whenever $x_n\to x$ and 
$y_n\to y$, then the $\phi$-contractivity of the members of the sequence $(T_n)$ implies the 
$\phi$-contractivity of $T_0$. Indeed, taking the limit $n\to\infty$ in the inequalities
$d(T_nx,T_ny)\le\phi\bigl(d(x,y)\bigr)$ and using the self-continuity of $d$, the inequality 
$d(T_0x,T_0y)\le\phi\bigl(d(x,y)\bigr)$ follows. Note that the self-continuity of $d$ holds automatically in 
metric spaces and also in ultrametric spaces. 

\section{Applications and Concluding Remarks}

Not claiming completeness, let us present here some immediate consequences of \thm{mainfix} and \thm{stabfix},
respectively. For their proof, one should observe only that in each case the underlying semimetric spaces
are regular. Moreover, the extra assumption of \thm{stabfix} is obviously satisfied. 

\Cor{Matkowski1}{If $(X,d)$ is a complete $c$-relaxed metric space or complete $c$-inframetric space and
$\phi$ is a comparison function, then every $\phi$-contraction has a unique fixed point.}

\Cor{Matkowski2}{Let $(X,d)$ be a complete metric space or a complete ultrametric space. If $(T_n)$ is a 
sequence of $\phi$-contractions converging pointwise to $T\colon X\to X$, then $T$ is a 
$\phi$-contraction, and the sequence of the fixed points of $(T_n)$ converges to the fixed point of $T$.}

Note, that the special case $c=1$ of \cor{Matkowski1} reduces to the result of Matkowski \cite{Mat75} which
is a generalization that of Browder \cite{Brow68}. Each of these results generalizes the Banach Fixed Point
Theorem under the particular choice $\phi(t)=qt$ where $q\in[0,1[$ is given. Let us also mention, that the
contraction principle was also discovered and applied independently (in a few years later after Banach) by
Caccioppoli \cite{Cac30}. Our last corollary demonstrates the efficiency of \thm{mainfix} in a particular
case when the Banach Fixed Point Theorem cannot be applied directly. For its proof, one should combine 
\cor{Matkowski1} with \lem{equiv}.

\Cor{extension}{If $(X,d)$ is a semimetric space with a semimetric $d$ equivalent to a complete $c$-relaxed 
metric or to a complete $c$-inframetric on $X$, then every $\phi$-contraction has a unique fixed point.}

As a final remark, let us quote here a result due to Jachymski, Matkowski, and \'Swi\k{a}tkowski, which
is a generalization of the Matkowski Fixed Point Theorem (for precise details, consult \cite{JacMatSwi95}).

\Thmn{Assume that $(X,d)$ is a complete Hausdorff semimetric space such that there is some $r>0$ for which the
diameters of balls with radius $r$ are uniformly bounded and in which the closure operator induced by $d$ is
idempotent. If $\phi$ is a comparison function, then every $\phi$-contraction has a unique fixed point.}

The proof of this theorem is based on the fundamental works of Chittenden \cite{Chi17} and Wilson \cite{Wil31}.
In view of \lem{regsem}, the first two conditions of the result above is always satisfied in regular semimetric
spaces. However, the exact connection between the properties of the basic triangle function and the idempotence
of the metric closure has not been clarified yet. The interested Readers can find further details on the topology
of semimetric spaces in the papers of Burke \cite{Bur72}, Galvin and Shore \cite{GalSho84}, and by McAuley
\cite{Mca56}.

%\bibliographystyle{amsplain}
%\bibliography{fixponthoz}
%\end{document}

\end{document}